\documentclass[12pt,a4paper]{article}
\usepackage{amsmath}
\usepackage{amsfonts}
\newtheorem{theorem}{Theorem}[section]
\newtheorem{lemma}{Lemma}[section]

\newtheorem{proposition}{Proposition}[section]

\title{\bf On perfect order subsets in finite groups}

\author{Nguyen Trong Tuan 
\\{\small\em High School for gifted}\\{\small\em VNU - HCM City}\\ {\small\em 145 Nguyen Chi Thanh str., Dist. 5, Ho Chi Minh City, Vietnam}\\ {\small\em e-mail: 
tuanhoangnhi@yahoo.com}\\ Bui Xuan Hai\\
{\small\em Faculty of Mathematics $\&$ Computer Science, University of Science}\\{\small\em VNU - HCM City}\\ {\small\em 227 Nguyen Van Cu str., Dist. 5, Ho Chi Minh City, Vietnam}\\ {\small\em e-mail: bxhai@hcmuns.edu.vn}}

\date{8 July 2010}
\rm\em 
\normalsize
\begin{document}
\maketitle
\newcommand{\dpcm}{ \hfill \rule{3mm}{3mm}}
\def\Box{\dpcm}
\def\xd{\linebreak} 
\begin{abstract} If $G$ is a finite group and $x\in G$ then the set of all elements of $G$ having the same order as $x$ is called {\em an order subset of $G$ determined by $x$} (see \cite{Car1}). We say that $G$ is a {\em group with perfect order subsets} or briefly, $G$ is a {\em $POS$-group} if the number of elements in each order subset of $G$ is a divisor of $|G|$. In this paper we prove that for any $n\geq 4$, the symmetric group $S_n$ is not $POS$-group. Together with the result in \cite{Das}, this gives the complete positive answer to Conjecture 5.2 in \cite{Car2}.
 
\end{abstract}



\newpage

\section{Introduction}

Throughout this paper, all considered groups are finite and for a  group $G$ we denote by  $|G|$  the order of  $G$, while for an element $x\in G$, the order of $x$ is denoted by $o(x)$. We  denote also by $\mathbb{N}$ and $\mathbb{Z}^+$ the sets of all non-negative and positive integers respectively. If $m\in\mathbb{Z}^+$ then $G^m$ denotes the direct product $G\times G\times\ldots\times G$. In a group $G$, define the following equivalence relation:
$$x\sim y\Longleftrightarrow o(x)=o(y).$$
The equivalence class defined by an element $x$ is denoted  by $\overline{x}$ and is called {\em an order subset} of $G$. Following the work \cite{Car1}, we say that $G$ is a {\em group with perfect order subsets} or briefly, $G$ is a {\em $POS$-group} if the number of elements in each order subset of $G$ is a divisor of $|G|$. In \cite{Car1}, the authors study properties of some abelian $POS$-groups and they established some curious  connection of such groups and Fermat numbers. In \cite{Car2}, the authors have extended their study for non-abelian groups and they have obtained some interesting properties for such groups. Also, in this work, some examples of non-abelian $POS$-groups are given. In particular, it is obvious that the symmetric group $S_3$ on three letters is a non-abelian $POS$-group. However, the authors conjectured that for $n\geq 4, A_n$ and $S_n$ do not have perfect order subsets, i.e. they are not $POS$-groups. Recently, in \cite{Das}, Ashish Kumar Das have proved this conjecture for groups $A_n$. Our main purpose in this paper is to prove that the conjecture for groups $S_n$ is also true. As an useful  additional information, in Section 2 we give some examples of groups having no perfect order subsets.

\section{Examples of groups having no perfect order subsets}

In this section we give some examples of groups not necessarily abelian, having no perfect order subsets. In the first, we note that the direct  product of $POS$-groups does not necessarily be a $POS$-group as the following proposition shows.

\begin{proposition}\label{pro2.1}  For $\alpha, t\in\mathbb{Z}^+, t> 1, (\mathbb{Z}_{2^{\alpha}})^t$ is not a $POS$-group.
\end{proposition} 

\noindent
{\bf Proof.} The order of any element of $(\mathbb{Z}_{2^{\alpha}})^t$ is of the form $2^i, i\leq \alpha t$. By [2, Lemma 1], the number of elements of the order $2^{\alpha}$ is $(2^{\alpha-1})^t(2^t-1)$. Since $t> 1, (2^t-1)$ does not divide $2^{\alpha t}$. Therefore $(\mathbb{Z}_{2^{\alpha}})^t$ is not a $POS$-group.
 \dpcm

For $2n\geq 4$, denote by $D_{2n}$ the dihedral group, defined by the following:
$$D_{2n}:=\langle a, b \vert a^n=1, b^2=1, bab^{-1}=a^{-1}\rangle.$$

\begin{lemma}\label{lem2.1} If $n$ is an even integer then $D_{2n}$ is not a $POS$-group.
\end{lemma}

\noindent
{\bf Proof.} In a dihedral group $D_{2n}:=\langle a, b \vert a^n=1, b^n=1, aba^{-1}=a^{-1}\rangle$, for every $i, 0\leq i< n$, we have $(a^ib)^2=1$. So $D_{2n}$ contains $n$ elements of order $2$ of this form. Now, if $n$ is even, then $a^{n/2}$ is  also an  element of order $2$. So the number of elements of order $2$ is not a divisor of $2n$ and it follows that $D_{2n}$ is not a $POS$-group.\dpcm

\begin{lemma}\label{lem2.2} If there exist at least two odd prime divisors of $n$ then $D_{2n}$ is not a $POS$-group.
\end{lemma}

\noindent
{\bf Proof.}  Suppose that $n=\prod_{k=1}^r p_k^{\alpha_k}, r\geq 2, p_k$ are all odd primes. Then, the number of elements of the order $n$ is $\varphi(n)=\prod_{k=1}^r p_k^{\alpha_k-1}(p_k-1)$. It follows that $\prod_{k=1}^r (p_k-1)$ is a divisor of $2\prod_{k=1}^r p_k$. But, by our assumption this is a contradiction.\dpcm

\begin{theorem}\label{th2.1} $D_{2n}$ is a $POS$-group if and only if $n=3^{\alpha}, \alpha\in \mathbb{Z}^+$.
\end{theorem}

\noindent
{\bf Proof.} Suppose that  $D_{2n}$ is a $POS$-group. In view of lemmas \ref{lem2.1} and \ref{lem2.2}, $n=p^\alpha$ for some odd prime $p$.  Then $b^p=b\neq 1$. It follows that every element of a order $p^\alpha$ is of the form $a^i$ with $(i, p)=1$ and $1\leq i< p^\alpha$. So, the number of elements of a order $p^\alpha$ is $\varphi(p^\alpha)=(p-1)p^{\alpha-1}$. Since this number is a divisor of $2p^\alpha$, it follows $p=3$. 

Conversely, suppose that $n=3^{\alpha}$. Then, the order of any element of $D_{2n}$ is of the form $2^i. 3^\beta$, where $\beta\leq \alpha$ and $i\in \{0, 1\}$. 

If $i=0$ then the number of elements of a order $3^\beta$ is $\varphi(3^\beta)=2.3^{\beta-1}$, which is a divisor of $n=2.3^\alpha$. Now, if $i=1$ then any element of a order $2.3^\beta$ must be of the form $a^kb, 0\leq k< n$. Then we have 
$$(a^kb)^2=a^kba^kb=a^kba^kb^{-1}=a^k(ba^kb^{-1})=a^ka^{-k}=1.$$ 
It follows that $\beta=0$. Thus, such elements have the order $2$ and there are exactly $n$ such elements. Hence $D_{2n}$ is a $POS$-group.\dpcm

Recall that for $n\geq 3$, the generalized quaternion group $Q_n$ is defined by the following:
$$Q_n:=\langle a, b \vert a^{2^{n-1}}=1, b^2=  a^{2^{n-2}}, bab^{-1}=a^{-1}\rangle.$$

Generalized quaternion groups are non-abelian non- $POS$ groups. In fact, we have the following result:

\begin{proposition}\label{pro2.2} For $n\geq 3$, a generalized quaternion group $Q_n$ is not a $POS$-group.
\end{proposition}

\noindent
{\bf Proof.} Consider a generalized quaternion group 
$$Q_n:=\langle a, b \vert a^{2^{n-1}}=1, b^2=  a^{2^{n-2}}, bab^{-1}=a^{-1}\rangle.$$
Since $bab^{-1}=a^{-1}$, $ba^ib^{-1}=a^{-i}$ for all $i, 0\leq i<2^{n-1}$. From the last equality it follows
$$(a^ib)^2=a^i(ba^ib^{-1})b^2=a^ia^{-i}b^2=b^2=a^{2^{n-2}}.$$
Hence $(a^ib)^4=(a^{2^{n-2}})^2=a^{2^{n-1}}=1$. So, there are $2^{n-1}$ elements of the order $4$ of the form $a^ib,  0\leq i<2^{n-1}$. On the other hand, the order of the element $a^{2^{n-3}}$ is $4$. Hence, the number of elements of order $4$ does not divide $2^n$. So $Q_n$ is not $POS$-group.\dpcm
 
\section{The symmetric groups}

As we have mentioned in the Introduction, the symmetric group $S_3$ is a $POS$-group. In [3, Conjecture 5.2], the authors conjectured that  for any $n\geq 4$, the group $S_n$ is not a $POS$-group. Our main purpose in this section is to  give the positive answer to this conjecture. In fact, we shall prove the following result:

\begin{theorem}\label{th3.1} For any integer $n\geq 4$, the symmetric group $S_n$ is not a $POS$-group.
\end{theorem}

To prove this theorem, we need some lemmas.

\begin{lemma}\label{lem3.1} Let $p$ be an odd prime number. If $n=2p+r$ with $r\in\{0, 1, \ldots, p-1\}$, then $S_n$ is not a $POS$-group.
\end{lemma}

\noindent
{\bf Proof.} Suppose that under our supposition, $S_n$ is a $POS$-group. Consider any element $\alpha$ of the order $p$ in $S_n$. Then, either $\alpha$ is a cycle of the length $p$ or it is a product of two disjoint cycles of the same length $p$. For the convenience, we call $\alpha$ an element of type $1$ for the first case and $\alpha$ an element of type $2$ for the second case. Obviously that the number of all elements of type $1$ is
 $$\frac{A^{p}_{n}}{p}=\frac{n!}{p(n-p)!}=\frac{n!}{p(p+r)!}$$
and the number of all elements of type $2$ is
$$\frac{1}{2} \times \frac{A^{p}_{n}}{p}\times \frac{A^{p}_{n-p}}{p}=\frac{1}{2} \times \frac{n!}{p(n-p)!}\times \frac{(n-p)!}{p(n-2p)!}=\frac{n!}{2p^{2}r!}.$$

Hence, the number of all elements of the order $p$ in $S_n$ is
$$d=\frac{n!}{p(p+r)!}+\frac{n!}{2p^{2}r!}=\frac{2pr!+(p+r)!}{2p^{2}r!(p+r)!}n!$$

Since $S_n$ is a $POS$-group, 
$$\frac {2p^{2}r!(p+r)!}{2pr!+(p+r)!}$$
is an integer or
$$\frac{2p^{2}(p+r)!}{2p+(r+1)(r+2)\ldots (p+r)}=\frac{2p^{3}r!(r+1)\ldots(p-1)(p+1)\ldots(p+r)}{p[2+(r+1)\ldots(p-1)(p+1)\ldots(p+r)]}$$
is an integer. Therefore 
$$\frac{2p^{2}r!(r+1)\ldots(p-1)(p+1)\ldots(p+r)}{2+(r+1)\ldots(p-1)(p+1)\ldots(p+r)}\eqno (1)$$
is an integer. Since
$$(r+1)\ldots(p-1)(p+1)\ldots(p+r) \equiv  (p-1)!\pmod p, $$
in view of Wilson's Theorem we have
$$ 2+ (r+1)\ldots(p-1)(p+1)\ldots(p+r) \equiv 1 \pmod p. \eqno (2)$$

From (1) and (2) it follows that 
$$\frac{2r!(r+1)\ldots(p-1)(p+1)\ldots(p+r)}{2+(r+1)\ldots(p-1)(p+1)\ldots(p+r)}=\frac{2r!A}{2+A}$$
is an integer, where $A=(r+1)\ldots(p-1)(p+1)\ldots(p+r)$. Since $gcd(A,2+A)=gcd(A,2)=1$ or $2, \displaystyle\frac{4r!}{2+A} $ is an integer, that is a contradiction in view of the following inequalities:
$$ 2+ A >(p+1)(p+2)(p+3)\ldots(p+r) > 4(p+2)\ldots(p+r) >4(1.2\ldots r)=4r!.$$

The proof is now complete.
\dpcm

\begin{lemma}\label{lem3.2} If $n =3p +r$, where $p$ is a odd prime and $r \in \{0,1,2,\ldots,p-1\}$, then $S_n$ is not a $POS$-group.
\end{lemma}

\noindent
{\bf Proof.}  Suppose that under our supposition, $S_n$ is a $POS$-group. Consider any element $\alpha$ of the order $p$ in $S_n$. Then, either $\alpha$ is a cycle of the length $p$ or it is a product of two or three disjoint cycles of the same length $p$. The number of elements in each of these cases is
$$\frac{A^{p}_{n}}{p}=\frac{n!}{p(n-p)!}=\frac{n!}{p(2p+r)!},$$
$$\frac{1}{2} \times \frac{A^{p}_{n}}{p}\times \frac{A^{p}_{n-p}}{p}=\frac{1}{2} \times \frac{n!}{p(n-p)!}\times \frac{(n-p)!}{p(n-2p)!}=\frac{n!}{2p^{2}(p+r)!}$$
and
$$\frac{1}{6} \times \frac{A^{p}_{n}}{p}\times \frac{A^{p}_{n-p}}{p}\times \frac{A^{p}_{n-2p}}{p}=\frac{1}{6} \times \frac{n!}{p(n-p)!}\times \frac{(n-p)!}{p(n-2p)!} \times \frac{(n-2p)!}{p(n-3p)!}=\frac{n!}{6p^{3}r!}$$
respectively. Hence, the number of elements of the order $p$ in $S_n$ is
 $$d=n!\left[\frac{1}{p(2p+r)!}+\frac{1}{2p^{2}(p+r)!}+\frac{1}{6p^{3}r!}\right].$$

Since $S_n$ is a $POS$-group, $d$ must be divided $n!$ and, consequently
$$k=\frac{6p^{3}.r!(p+r)!(2p+r)!}{6p^{2}(p+r)!r!+3p(2p+r)!r!+(2p+r)!(p+r)!}$$ is an integer. By setting $A:=(r+p+1)\ldots (2p-1)(2p+1)\ldots (2p+r)$ and the direct calculation we have
$$k=\frac{6p^{3}.(p-1)!(p+1)\ldots (p+r)A}{3+3A+(r+1)\ldots (p-1)(p+1)\ldots (p+r)A}.\eqno (3)$$

By applying of Wilson's Theorem we get
$$(r+1)\ldots (p-1)(p+1)\ldots (p+r) \equiv -1 \pmod p\eqno (4)$$
and, consequently we have
$$A=(p+r+1)\ldots (2p-1)(2p+1)\ldots (2p+r) \equiv -1 \pmod p.\eqno (5)$$

From (4) and (5) it follows that
$$3+3A+(r+1)\ldots (p-1)(p+1)\ldots (p+r)A \equiv 1 \pmod p .$$

Since 
$$gcd \left(A,3+A[3+(r+1)\ldots(p-1)(p+1)\ldots(p+r)]\right)=gcd(3,A)$$
and $k$ is an integer, it follows from (3) that
$$B:=\frac{18(p-1)!(p+1)\ldots (p+r)}{3 +3A+(r+1)\ldots (p-1)(p+1)\ldots (p+r)A}$$ 
is an integer. Now, we claim that 
$$(r+1)\ldots(p-1)(p+1)\ldots(p+r)A > 18(p-1)!(p+1)\ldots (p+r).$$

In fact, this inequality is equivalent to the following one:
$$A=(r+p+1)\ldots(2p-1)(2p+1)\ldots(2p+r)>18r!.$$

Since $p$ is an odd prime, the last inequality holds as the following calculation shows:

$A=(2p+1)(2p+2)\ldots(2p+r)(r+p+1)\ldots(2p-1)\geq (2p+1)(p+1)(2.3\ldots r)$

\hspace*{0.45cm}$=(2p+1)(p+1).r!>18r!.$

Clearly, what we have claimed shows that $B$ is not an integer. This contradiction completes the proof of the lemma.\dpcm

\begin{lemma}\label{lem3.3} If $n =4p $, where $p$ is a odd prime,  then $S_n$ is not a $POS$-group.
\end{lemma}

\noindent
{\bf Proof.}  Suppose that under our supposition, $S_n$ is a $POS$-group. Let $d$ be the number of elements of the order $p$ in $S_n$. Then we have

$d=\displaystyle\frac{A^{p}_{n}}{p}+\frac{1}{2}\times  \frac{A^{p}_{n}}{p} \times\frac{A^{p}_{n-p}}{p}+\frac{1}{6} \times \frac{A^{p}_{n}}{p}\times \frac{A^{p}_{n-p}}{p} \times\frac{A^{p}_{n-2p}}{p}+\frac{1}{24} \times \frac{A^{p}_{n}}{p} \times\frac{A^{p}_{n-p}}{p} \times\frac{A^{p}_{n-2p}}{p}\times \frac{A^{p}_{n-3p}}{p}$

\hspace*{0.3cm}$=n!\displaystyle\left[\frac{1}{p(3p)!}+\frac{1}{2p^{2}(2p)!}+\frac{1}{6p^{3}p!}+\frac{1}{24p^{4}}\right]$

\hspace*{0.3cm}$=\displaystyle \frac{24p^{3}.(2p)!p!+12p^{2}.(3p)!p!+4p.(3p)!(2p)!+p!(2p)!(3p)!}{24p^{4}.p!(2p)!(3p)!}.n!$

\hspace*{0.3cm}$=\displaystyle\frac{24p^{3}+12p^{2}(2p+1)\ldots (3p)+4p(p+1)\ldots (3p)+(3p)!}{24p^{4}(3p)!}.n!.$

Since $d$ divides $n!$, 
$$\displaystyle\frac{n!}{d}=\frac{24p^{4}(3p)!}{24p^{3}+12p^{2}(2p+1)\ldots (3p)+4p(p+1)\ldots (3p)+(3p)!}$$
is an integer.  By dividing both numerator and denominator of the last fraction by $6p^3$ we get
$$\frac{n!}{d}=\frac{24p^{4}.(p-1)!(p+1)\ldots(2p-1)(2p+1)\ldots(3p-1)}{4+(2p+1)\ldots(3p-1)[6+4(p+1)\ldots(2p-1)+(p-1)!(p+1)\ldots(2p-1)]}.$$ 

By setting  $$A=(2p+1)\ldots(3p-1)$$ and  $$M=4+A[6+4(p+1)\ldots(2p-1)+(p-1)!(p+1)\ldots(2p-1)],$$  we have
 $$\frac{n!}{d}=\frac{24p^{4}(p-1)!(p+1)\ldots(2p-1)A}{M}.\eqno (6)$$

In view of Wilson's Theorem we have
$$(2p+1)\ldots(3p-1)=(2p+1)(2p+2)\ldots(2p+p-1) \equiv (p-1)! \equiv -1 \pmod p;$$
$$ (p+1)\ldots(2p-1)=(p+1)(p+2)\ldots(p+p-1)\equiv (p-1)! \equiv -1 \pmod p.$$

Hence $M \equiv 4+(-1)[6-4+(-1)(-1)] \equiv 1 \pmod p$. Consequently, $gcd(M,p^4)=1$ . Therefore, in view of (6) we conclude that 
$$\frac{24(p-1)!(p+1)\ldots(2p-1)A}{M}$$
is an integer. Note that, if $m$ is a common divisor of $A$ and $M$, then $m$ must divide $4$. In particular, $gcd(A, M)$ must be $1, 2$ or $4$. It follows that
 $$C:=\frac{96(p-1)!(p+1)\ldots(2p-1)}{M}$$
is an integer. However, we can check that $C$ is not an integer for any odd prime $p$. In fact, if $p=3$, then
$$C=\frac{3840}{7060}$$
which is not an integer. Now, suppose that $p\geq 5$. Then we have
$$A=(2p+1)\ldots(3p-1) \geq (2.5+1).12.13.14 > 96$$
and
$$M > (p-1)!(p+1)\ldots(2p-1)A > 96(p-1)!(p+1)\ldots(2p-1).$$

Hence, in this case $C$ is not an integer too. This  contradiction completes the proof of the lemma. \dpcm

Now, we are ready to prove the main theorem in this section.\\

\noindent{\bf Proof of Theorem 3.1.}

For $n=6$ and $n=7$, the  desired result follows from Lemma \ref{lem3.1} by taking $p=3, r=0$ and $p=3, r=1$ respectively.  For $n=4$ and $n=5$, note that the number of elements of the order $2$ in $S_4$ and $S_5$ is $9$ and $25$ respectively. So, $S_4$ and $S_5$ are both non-$POS$ groups.

Now, suppose that $n\geq 8$ and $\displaystyle m=\left[\frac{n}{4}\right]$. According to Bertrand's Postulate (see, for example [4, Theorem 5.8, p. 109]), there exists some  prime $p$ such that 
$$ m < p < 2m .$$ 

Note that
$$p <2m =2\left[\frac{n}{4}\right]\leq 2\frac{n}{4}=\frac{n}{2}.$$

If $\displaystyle \left[\frac{n}{4}\right] < p \leq \frac{n}{4}$,  then $n=4p$ and the conclusion follows from Lemma 3.3. Therefore, we can suppose that 
$$\displaystyle \frac{n}{4} <p < \frac{n}{2}.$$

If $\displaystyle \frac{n}{4} < p \leq \frac{n}{3}$, then  $n =3p+r$ with $r \in \{0,1,2,\ldots,p-1\}$ and the conclusion follows from Lemma \ref{lem3.2}. If
$\displaystyle \frac{n}{3} < p < \frac{n}{2},$ then $n=2p+r$ with $r\in \{0,1,2,\ldots,p-1\}$ and the conclusion follows from Lemma \ref{lem3.1}. The proof of the theorem is now complete. \dpcm

\end{document}